\documentclass[12pt]{amsart}
\usepackage{mathrsfs}
\usepackage{amsmath}
\usepackage{verbatim}
\usepackage{hyperref}
\usepackage{amsfonts}
\usepackage{bbm}
\usepackage{amscd}
\usepackage{color}
\usepackage{amssymb}
\usepackage[all, knot]{xy}
\usepackage[x-cp1252]{inputenx}
\usepackage{CJK}
\xyoption{arc} 

\usepackage[foot]{amsaddr}
\usepackage{amstext}
\usepackage{amsxtra}
\usepackage{amsmath}
\usepackage{amsopn}
\usepackage{amsthm}
\newtheorem{thm}{Theorem}[section]
\theoremstyle{definition}
\newtheorem{defn}{Definition}[section]
\theoremstyle{plain}
\newtheorem*{mr}{Main Result}
\theoremstyle{remark}

\theoremstyle{plain}
\newtheorem{lem}[thm]{Lemma}
\theoremstyle{plain}
\newtheorem{cor}[thm]{Corollary}
\theoremstyle{plain}
\newtheorem{prop}[thm]{Proposition}
\newtheorem{con}{Conjecture}[section]
\newcommand{\mc}{\mathcal}

\newcommand{\e}{\ensuremath{\mathcal{E}^{\vee}}}    
\newcommand{\Sym}{\rm{Sym}}  
\newcommand{\s}{\rm{Sym}}  
\CompileMatrices
\title{A Correlator Formula for Quantum Sheaf Cohomology}
\author{Zhentao Lu}
\address{Mathematical Institute, University of Oxford, 
Oxford OX2 6GG, UK}
\email{Zhentao.Lu@maths.ox.ac.uk}
\begin{document}
\begin{abstract}
For a class of monadic deformations of the tangent bundles over nef-Fano smooth projective toric varieties, we study the correlators using quantum sheaf cohomology. We prove a summation formula for the correlators, confirming a conjecture by McOrist and Melnikov in physics literature. This generalizes the Szenes-Vergne proof of Toric Residue Mirror Conjecture for hypersurfaces.


\end{abstract}
\maketitle
\section{Introduction}

The study of quantum sheaf cohomology (QSC) arises from the physics problem of understanding the Gauged Linear Sigma Model (GLSM) studied by Witten \cite{witten1993phases}. There are two versions of GLSM, the $(2, 2)$ theory and the $(0, 2)$ theory, where $(2, 2)$ and $(0, 2)$ indicates the amount of supersymmetry of the theory. Both theories have different phases. Mathematically speaking, in geometric phases, both theories study maps from Riemann surfaces to compact K\"ahler manifolds. The former theory considers the manifold with its tangent bundle, while the latter considers more general vector bundles. QSC comes out of the understanding of $(0, 2)$ GLSMs in this geometric setting.

The study of the $(2, 2)$ GLSM is more mature and the associated quantum cohomology theory is studied by Batyrev \cite{batyrev1993quantum} and Morrison-Plesser \cite{morrison1995summing}. The main result there, the Toric Residue Mirror Conjecture (TRMC, see Equation (\ref{trmc})), is formulated in \cite{batyrev2002toric} based on \cite{morrison1995summing}, and is proved independently by Szenes-Vergne \cite{szenes2004toric} and Borisov \cite{borisov2005higher}. 

Quantum sheaf cohomology arose from the study of the $(0, 2)$ theory, which is relatively new and many problems remain open. See \cite{melnikov2012recent}\cite{mcorist2011revival} for surveys. Mathematically, the basic object studied here is a compact K\"ahler manifold $V$ with an omalous\footnote{ ``Omalous" means ``non-anomalous", i.e. the Chern classes satisfy $c_i(\mathcal{E}) = c_i (V), i=1,2$.} holomorphic vector bundle $\mathcal{E}$ on it. QSC is then a `quantum deformation' of the classical sheaf cohomology ring $H^*(V, \wedge^*\mc{E})$.  While a general QSC theory is still missing, Donagi, Guffin, Katz and Sharpe in \cite{donagi2014mathematical} defined QSC for any smooth projective toric variety $X$ together with a bundle $\mathcal{E}$ defined by the deformed toric Euler sequence (\ref{euler1}), based on previous work \cite{katz2006notes}\cite{guffin2010deformed}. We will recall this in Section \ref{quantum-correlator}. Bundles defined this way are naturally omalous, and they can be studied using Koszul complex. The quantum sheaf cohomology ring is defined by specifying the quantum Stanley-Reisner ideals. This enables the authors of \cite{donagi2014mathematical} to define the quantum correlators with values in a one-dimensional complex vector space $H^*$.

An important object in the $(0, 2)$ theory (as well as in the $(2, 2)$ theory) is the set of \textit{correlators}. For cohomology elements in $H^1(V,\mathcal{E}^\vee)$, the classical correlator is the cup product, which can be interpreted as a sheaf cohomology analog of the intersection number of divisors, while the quantum correlator is a weighted sum of classical correlators of the moduli spaces parametrized by effective curve classes of $V$.


From the physics side, McOrist and Melnikov formulated a conjecture about the quantum correlators in \cite{mcorist2008half}. 
\begin{con}\label{toric_case}
For a toric variety $V$ with an (automatically omalous) holomorphic vector bundle $\mathcal{E}$ defined by a deformed toric Euler sequence (See (\ref{euler1}) below), the quantum correlator of $\sigma_i$'s in $H^1(\mathcal{E}^\vee)$ can be computed by the following summation formula:

\begin{equation}\label{RHS0}
\langle \sigma_{i_1},...,\sigma_{i_s}\rangle^{quantum} =\displaystyle \sum_{\{u\in W^\vee|\tilde{v}_j(u)=q_j\}} \frac{\sigma_I}{\prod_{c\in [\Sigma(1)]} Q_c} \frac{\prod_{j=1}^{r}\tilde{v}_j}{\det_{j,k} (\tilde{v}_{j,k})}.
\end{equation}
\end{con}
In the above formula, the quantities on the right hand side are constructed from the map in the deformed toric Euler sequence defining $\mathcal{E}$. We will give the precise definition of the notations in Section \ref{quantum-sum}.

The authors of \cite{mcorist2008half} work in the physics theory of Coulomb branch and derive the formula (\ref{RHS0}) there. They then conjecture that the same formula holds for the geometric case as described in Conjecture \ref{toric_case}.

Conjecture \ref{toric_case} has the feature that the quantum correlators take values in the complex numbers. Compared to the result in \cite{donagi2014mathematical}, it offers a specific trivialization of the top cohomology, and an explicit way to compute it.

In this paper we prove Conjecture \ref{toric_case} for nef-Fano smooth projective toric varieties. The outline of the paper is as follows:


In section \ref{preliminaries} we set up the basic notations, introduce the classical correlators, and review some results in the case where $\mathcal{E}$ is the tangent bundle of the toric variety.

In Section \ref{classical}, we prove an integral formula that computes the classical correlators when $\mathcal{E}$ is a deformation of the tangent bundle $T_V$ with small deformation parameters.

In Section \ref{quantum-integral},  we define the quantum correlators following \cite{donagi2014mathematical}. Then we use this to rewrite the quantum correlator, which is the sum of all contributions of classical correlators from different moduli spaces parametrized by effective curve classes of the toric variety, into an integral form.

Section \ref{QSC-system} generalizes the corresponding result in \cite{donagi2015global-aspects}. In this section we study the number of solutions to the QSC relations. The result is used to guarantee that we have a finite sum in the summation formula, and to prove the properness lemma, which is an important ingredient of the proof of the main result.

In Section \ref{quantum-sum}, we prove Conjecture \ref{toric_case} by showing that the right hand side of Equation (\ref{RHS0}) can be written as an integral, which has the same integrand as the one in Section \ref{quantum-integral}. However the two integrals are over different cycles. We will then prove the equality of the two integrals by 
showing that the integration cycles contains homologous parts and extra parts where the integration vanishes. The proof is inspired by the $(2, 2)$ case work \cite{szenes2004toric} proving the TRMC. In the end we comment on the significance of this summation formula.

The main results of this paper were obtained in the author's Ph.D. thesis  
, but the presentation was improved based on subsequent ideas developed in \cite{donagi2015global-aspects}.

\textbf{Acknowledgements:} I would like to thank my thesis advisor, Ron Donagi, for his numerous  suggestions in all stages of writing this paper and for his constant encouragement. I also thank Xenia de la Ossa, Sheldon Katz and Ilarion Melnikov for helpful discussions. I acknowledge the  partial support of the Benjamin Franklin Fellowship and the Dissertation Completion Fellowship of the University of Pennsylvania, and the EPSRC grant EP/J010790/1.

\section{Preliminaries}\label{preliminaries}
\subsection{The basic setting}
Let $V$ be an $n$-dimensional smooth projective toric variety over $\mathbb{C}$, with the toric fan $\Sigma$. $\Sigma(k)$ is the collection of all $k$-dimensional cones of $\Sigma$. Each ray $i \in \Sigma(1)$ corresponds to a prime divisor $D_i$ via the orbit-cone correspondence. Let $M$ be the character lattice of $V$ and $N$ be the dual lattice. $v_i\in N_\mathbb{R}$ denotes the first integral point of the ray $i$. For a detailed treatment of toric varieties, see \cite{cox2011toric}.

Recall that the tangent bundle $T_V$ of $V$ fits in the toric Euler sequence
\begin{equation}
0 \to \mathcal{O}\otimes W^\vee \xrightarrow{E_0^\vee} \displaystyle\bigoplus_{i\in\Sigma(1)} \mathcal{O}(D_i) \to T_V \to 0,
\end{equation}
where $W$ is an $r$ dimensional complex vector space such that $\mc{O}\otimes W^\vee \cong \mc{O}^r$ (we choose $W^\vee$, the dual of $W$ to simplify the notations). By a \textit{deformed tangent bundle} we mean a holomorphic vector bundle
$\mathcal{E}$ on $V$,  defined by the deformed toric Euler sequence (sometimes referred as a \textit{monadic deformation}, since $\mathcal{E}$ is the cokernel of two bundles that are direct sums of line bundles)

\begin{equation}\label{euler1}
0 \to \mathcal{O}\otimes W^\vee \xrightarrow{E^\vee} \displaystyle\bigoplus_{i\in\Sigma(1)} \mathcal{O}(D_i) \to \mathcal{E} \to 0.
\end{equation}

And we will mainly make use of the dual sequence

\begin{equation}\label{euler2}
0 \to \e \to \displaystyle\bigoplus_{i\in\Sigma(1)} \mathcal{O}(-D_i) \xrightarrow{E} \mathcal{O}\otimes W \to 0.
\end{equation}

\subsection{The definition of the polymology}

Note that $W\cong H^1(\mathcal{E}^\vee)$ for all deformed tangent bundles $\mathcal{E}^\vee$.
The sum of sheaf cohomologies $\bigoplus_{p,q} H^q(V, \wedge^q\e)$ together with cup product forms an associative algebra $H_\mathcal{E}^*(V)$ called the \textit{polymology} of $\mathcal{E}$ (\cite{donagi2014mathematical}, Definition 1.1).

Let $SR(V)$ be the \textit{Stanley-Reisner} ideal of $V$. $SR(V) =\langle x_{I}| I\subset \Sigma(1)$ and the rays in $I$ do not form a cone $ \rangle$.
Recall that the cohomology of toric varieties can be described using Stanley-Reisner ideals:

\begin{equation}
H_{T_V}^*(V)\cong {\rm Sym}^*W / SR(V).
\end{equation}
We can describe the polymology in a similar way.

Let $E$ be the second map in (\ref{euler2}) whose kernel defines $\mathcal{E}^\vee$. Then $E$ is in $Hom(\oplus(-D_i), \mathcal{O}\otimes W) \cong H^0(\oplus\mathcal{O}(D_i))\otimes W$. Let $E =\sum_i E_i$ and $\Delta_i = \{m\in M_\mathbb{R} | \langle m, v_j\rangle \geq -\delta^i_j, $ for all $j\in \Sigma(1) \}$. For each $i\in \Sigma(1)$ we have an expression in monomials
\begin{equation}
E_i = \sum _{m\in \Delta_i\cap M} a_{im}\chi^m,
\end{equation}
where $\chi^m, m\in \Delta_i\cap M$ is a basis of $H^0(\mathcal{O}(D_i))$ and $a_{im}$ takes values in $W$.

Let $S(i)=\{i^\prime\in\Sigma(1) \ |\   D_{i^\prime}\sim D_i\}$.
Since it is shown in \cite{donagi2014mathematical} that the quantum sheaf cohomology does not depend on non-linear deformations, we can focus on the linear part of $E_i$ defined by
\begin{equation}
E_i^{\rm lin} := \sum_{i^\prime\in   S(i)   } a_{ii^\prime}x_{i^\prime},
\end{equation}
where $D_{i^\prime}\sim D_i$ means they are linearly equivalent divisors, and $x_{i^\prime} \in H^0(\mathcal{O}(D_i))$ is the element in the homogeneous coordinate ring of the toric variety corresponding to a global section vanishing on $D_{i^\prime}$.

For the divisor class $c$, we introduce the notation 
\begin{equation}\label{Qc}
 Q_c = \det (a_{ij})
 \end{equation} with $i,j$ running through all rays $k\in\Sigma(1)$ such that $[D_k] = c$.

An important notion in toric geometry is the primitive collection. A primitive collection $P\subset \Sigma(1)$ is a collection of rays such that no cone in $\Sigma$ contains all the rays in $P$, but for any proper subset $P^\prime$ of $P$,  there is a cone in $\Sigma$ that contains all the rays in $P^\prime$.

It can be shown that for each primitive collection $P$ of rays in $\Sigma(1)$, if it contains one ray $i$, it has to contain $S(i)$ (\cite{donagi2014mathematical}). 

For any subset $S\subset \Sigma(1)$, let $[S]=\{[D_i]\ |\ i\in S\}$. 


Define the deformed Stanley-Reisner ideal $SR(V,\mathcal{E}) \subset {\rm Sym}^*W$ to be the ideal generated by $\prod_{c\in [P]} Q_c$ with $P$ running through all primitive collections of the fan, i.e.,
\begin{equation}\label{sr}
SR(V,\mathcal{E}) = \langle\prod_{c\in [P]} Q_c \ |\ {P {\rm \ is \ a \  primitive\ collection}}\rangle.
\end{equation}

Then it is proved in \cite{donagi2014mathematical} that the polymology of $\mathcal{E}$ satisfies:

\begin{equation}\label{polymology}
H_\mathcal{E}^*(V)\cong {\rm Sym}^*W / SR(V,\mathcal{E}).
\end{equation}
\subsection{The classical correlators}
For $i\in \Sigma(1)$, let $\alpha_i\in W$ be the first Chern class of the toric invariant divisor $D_i$ (under the identification $W\cong H^1(\Omega)$) and denote $\mathfrak{U}=\{\alpha_i\ |\  i\in \Sigma(1)\}$. Let $\sigma_i \in W$ be a general element of $W$ and $\sigma_I =\prod_{i\in I}\sigma_i \in {\rm Sym}^*W$.

Note that (\ref{euler2}) implies that $\wedge^n \e \cong \mathcal{O}(-\sum_{i\in\Sigma(1)}D_i)\cong K$, the canonical bundle. Hence $H^n(V, \wedge^n\e)$ is one-dimensional. Identify $H^n(V, \wedge^n\e)$ with $\mathbb{C}$ by integrating over the fundamental class.
For $\sigma_i, i\in I$, one can first take the image of $\sigma_I$ in $H^*_{\mathcal{E}}(V)$, project to the degree $(n,n)$ part $[\sigma_I]_n$, and 
define the \textit{(classical) correlator} of $\sigma_i, i\in I$ to be the image of $[\sigma_I]_n$ in $\mathbb{C}$. Denote the correlator of $\sigma_i, i\in I$ by $\langle\sigma_I\rangle$.

\subsection{An integral formula for the $(2, 2)$ classical correlators}\label{$(2, 2)$ classical correlators}
Following the physicists' language, we call the case in which $\mathcal{E}
$ is the tangent bundle $T_V$ \textit{the $(2, 2)$ case}. In this case the polymology of $\mathcal{E}$ is just the cohomology of holomorphic forms, and the ring structure can be computed by intersection theory (recall that $V$ is a smooth projective toric variety).

In this section, we present an integral formula for the $(2, 2)$ classical correlators found by Szenes and Vergne in \cite{szenes2004toric}.
We will generalize this formula in Section \ref{classical}.

First we need more notations.

Choose a maximal cone $\sigma\in\Sigma(n)$, and fix an order of  $\alpha_{i_1},...,\alpha_{i_r}$ corresponding to rays that are not in $\sigma$. This fixes a translation invariant measure $d\mu$ on $W^\vee$, where $d\mu=d\alpha_{i_1}\wedge d\alpha_{i_2}\wedge...\wedge d\alpha_{i_r}$.

For each prime toric divisor class $c\in [\Sigma(1)]$, let $H_c$ be a hypersurface in $W^\vee$ defined by  $H_c=\{u\in W^\vee\ |\ Q_c(u)=0\}$.
Let $U(\mathcal{E})$ be the complement of the union of $H_c, c\in [\Sigma(1)]$. In the current case where $\mathcal{E}=T_V$, $U(T_V)$ is the complement of the union of hyperplanes defined by $\alpha_i=0, i\in \Sigma(1)$. Let $r=\dim W^\vee$.


We can then state the following theorem:

\begin{thm}\label{$(2, 2)$classical}There is a homology class  $h(T_V)\subset H_r(U(T_V),\mathbb{Z})$ such that the following integral computes the $(2, 2)$ classical correlators for any $\sigma_i\in W, i\in I$.
\begin{equation}
\langle \sigma_I\rangle =\frac{1}{(2\pi i)^r} \int_{h(T_V)} \frac{\sigma_I}{\prod_{i\in \Sigma(1)} \alpha_i}d\mu.
\end{equation}

Moreover, the homology class is represented by a disjoint union of tori with orientations, as described below.
\end{thm}
\subsubsection*{Description of the homology class.}
To describe $h(T_V)$, we first introduce the set
$\mathcal{FL}(\mathfrak{U})$ of complete flags 
\begin{equation}
F=\{F_0=\{0\}\subset F_1\subset F_2\subset ... \subset F_{r-1}\subset F_r =W\}, 
\end{equation}
such that each $F_j$ is generated by some $\alpha_i\in \mathfrak{U}$.

We say an ordered basis $\gamma^F=(\gamma^F_1,...,\gamma^F_r)$ of $W$ is \textit{compatible} with $F$, if the following conditions hold:

(a) $\gamma_j^F $ is rationally generated by $\alpha_i \in \mathfrak{U}$,

(b) $\{\gamma^F_m\}_{m=1}^j$ is a basis of $F_j$,

(c) $d\gamma_1^F \wedge...\wedge d\gamma_r^F = d\mu$.

Fix a $\xi$ in the K\"ahler cone $\mathfrak{c}$ of $V$.
Let $\mathcal{FL}^{+}(\mathfrak{U},\xi)$ be the set of those flags $F=\{F_j\}\in \mathcal{FL}(\mathfrak{U})$ such that $\xi$ is in the interior of the cone  spanned by $\kappa_j$, $j=1,2,...,r$, where $\kappa_j = \sum_{\{\alpha_i\in\mathfrak{U}|\alpha_i\in F_j\}} \alpha_i$.

For each flag $F$, we always fix a compatible basis $\gamma^F$. Let $u_j = \gamma_j^F(u)$ for $u\in W^\vee$. 

Consider the torus $T_F(\mathbf{\epsilon}) =\{u\in W^\vee \ |\  |u_j| = \epsilon_j, j = 1,...,r\}$. Let 
\begin{equation}\label{Z(e)}
Z(\mathbf{\epsilon}) = \sum_{F\in \mathcal{FL}^{+}(\mathfrak{U},\xi)} \nu(F) T_F(\mathbf{\epsilon}),
\end{equation} where $\nu(F) =\pm 1$ depending on the orientation of $\kappa_j$. Szenes and Vergne prove the theorem by showing that $Z(\mathbf{\epsilon})$ represents $h(T_V)$ for $\mathbf{\epsilon}$ in some specific neighbourhood of $0$. We will specify the constraint for $\mathbf{\epsilon}$ in Section \ref{classical}.

\section{An integral formula for $(0, 2)$ classical correlators} \label{classical}
In this section we prove an integral formula which computes the classical correlators for $\mathcal{E}$ being a deformation of the tangent bundle $T_V$ with small deformation parameters. The statement is Theorem \ref{$(0, 2)$classical} below.
\subsection{The integral formula} \label{classical_correlators}
We fix the same translation invariant measure $d\mu$ on $W^\vee$ as in Section \ref{$(2, 2)$ classical correlators}. Recall that $U(\mathcal{E})$ is the complement of the union of all the hypersurfaces $H_c=\{u\in W^\vee\ |\ Q_c(u)=0\}$ in $W^\vee$, for $c\in [\Sigma(1)]$. $r=\dim W^\vee$.
Our first result generalizes the formula for $(2, 2)$ classical correlators:
\begin{thm}\label{$(0, 2)$classical}
For a small deformation $\mathcal{E}$ of the tangent bundle $T_V$,
there is a homology class  $h(\mathcal{E})\subset H_r(U(\mathcal{E}),\mathbb{Z})$ such that the following integral computes the $(0, 2)$ classical correlators for any $\sigma_i\in W, i\in I$:
\begin{equation}
\langle \sigma_I\rangle =\frac{1}{(2\pi i)^r} \int_{h(\mathcal{E})} \frac{\sigma_I}{\prod_{c\in [\Sigma(1)]} Q_c}d\mu,
\end{equation}
where $Q_c$'s are the factors of $SR(V,\mathcal{E})$ introduced in (\ref{Qc}).
Moreover, the cycle of integration $h(\mathcal{E})$ is represented by $Z(\epsilon)$ as described in (\ref{Z(e)}).
\end{thm}
\textbf{Remark:}

(a) There are constrains on $\epsilon_j$.That is $N \epsilon_i < \epsilon_{i+1}, i=1,2,...,r-1$, for a sufficiently large $N$, namely, $N$ is larger than $N_0(F)$ which could be chosen as follows: For all $l\in \Sigma(1)$, write $\alpha_l =\sum a_{li} \gamma^F_i$. Define $ N_0(F) = r\cdot (\max_{l} (\frac{1}{|a_{li}|})\cdot (\max_{i} |a_{li}|))$.

(b) The integral vanishes when $\sigma_I\in SR(V,\mathcal{E})$. This constitutes the main part of the proof.

(c) It is shown in \cite{donagi2014mathematical} that the quantum sheaf relations do not depend on the non-linear deformations. Hence the correlation functions, being linear functions vanishing on the ideal generated by the quantum sheaf relations, do not depend on the non-linear deformations.

\subsection{Preparatory results}

We state and prove some lemmas before we prove Theorem \ref{$(0, 2)$classical} in next section.
\begin{lem}\label{k}
Fix a flag $F$ and a collection $\{\alpha_l,l\in L\}$. If there exist a $k$, $1\leq k\leq r$ such that $F_k$ is not generated by elements of $\{\alpha_l,l\in L\}$ but every $F_i, i< k$ is, then $\prod_{l\in L}\alpha_L$ is non-zero on the region $\Omega=\{u;|u_k|\leq\epsilon_k, |u_i|=\epsilon_i,{\rm for \ } i\neq k \}$. 
\end{lem}

Proof:

For any $l\in L$, express $\alpha_l$ in terms of the basis $\gamma^F_i$, we have $\alpha_l =\sum a_{li}\gamma^F_i$. Order the sum in the descending order regarding $i$-index, and call the largest index $s$. Then $s\neq k$ by the definition of $k$.

Then, for $u\in \Omega, l\in L$,
\begin{equation}
\begin{array}{rl}
|\alpha_l(u)|&=|\sum_{i=1}^s a_{li}u_i|\\
&\geq |a_s|\cdot |u_s| -\sum^{s-1}_{i=1} |a_i|\cdot |u_i|\\
&\geq |a_s|\epsilon_s - \sum^{s-1}_{i=1} |a_i|\cdot \epsilon_i\\
&> \epsilon_{s-1}(N_0(F)\cdot |a_s| - r\cdot \max_i (|a_i|))\\
&>0,
\end{array}
\end{equation}
by the definition of $N_0(F)$. $\qed$

\begin{lem}\label{degenerate} Let $P\subset\Sigma(1)$ be a primitive collection.

(a)

When $\{\alpha_i, i\in L=\Sigma(1)-P-J\}$ does not form a basis compatible with the flag $F$,
\begin{equation*}
\int_{T_F(\epsilon)} \frac{\alpha_J}{\prod_{i\not\in P}\alpha_i}d\mu =0.
\end{equation*}
In particular, when $\{\alpha_i,i\in L\}$ does not generate $W$, we have
\begin{equation*}
\int_{T_F(\epsilon)}\frac{f}{\prod_{i\in L}\alpha_i}d\mu=0
\end{equation*}
for any $T_F(\epsilon)$.

(b)
When $\alpha_i, i\in \Sigma(1)-P-J$ form a basis compatible with both the flag $F$ and $F^\prime$,
$\int_{T_{F^\prime}(\epsilon)} \frac{\alpha_J}{\prod_{c\not\in P}Q_c} =sgn(F^\prime, F)\cdot\int_{T_{F}(\epsilon)} \frac{\alpha_J}{\prod_{c\not\in P}Q_c}$, where $sgn(F^\prime, F) =\pm 1$ depending on the orientation.
\end{lem}

Proof: 

(a) Let $L=\Sigma(1)-P-J$. $\{\alpha_i,i\in L\}$ satisfies the assumption of Lemma \ref{k}, hence $\prod_{l\in L}\alpha_L$, is non-zero on $\Omega$.

 Hence the integrand is defined on the region $\Omega$. Note that $T_F(\epsilon)=\partial \Omega$. Hence $\int_{T_F(\epsilon)}\frac{f}{\prod_{i\not\in P}\alpha_i} d\mu=\int_\Omega d(\frac{f}{\prod_{i\not\in P}\alpha_i}d\mu) =0$.

(b) In this case we will rename the indices of $\{\alpha_i, i\in\Sigma(1)-P-J\}$, so that $F_j$ is generated by $\{\alpha_1,...,\alpha_j\}$. If we write $\alpha_j(u) =\sum_{i=1}^{j}a_{ji}u_i$, then $|\alpha_j(u)/a_{jj}| = |u_j + \sum_{i=1}^{i=j-1} \frac{a_{ji}}{a_{jj}}u_i|$. This allows us to use a linear homotopy map to show that $|\alpha_j(u)/a_{jj}|=\epsilon_j$ is homotopic to $|u_j|=\epsilon_j$. Thus we conclude that $T_F(\epsilon) = \pm [T_\alpha]$, where $T_\alpha =\{u; |\alpha_j(u)/a_{jj}| =\epsilon_j\}$ and $[T_\alpha]$ is its homology class.

Recall that $Q_c$ is a deformation of $\prod_{j,[D_j]=c} \alpha_j$. We may group the $\alpha_j$'s for $j\in J$ in a similar fasion. Namely, we have 
\begin{equation}
\alpha_J = \prod_c (\prod_{j\in J,[D_j]=c}\alpha_j).
\end{equation}

 Since $J\subset \Sigma(1)-P$, we have \begin{equation}
\frac{\alpha_J}{\prod_{c\not\in P}Q_c} = \frac{1}{\prod_{c\not\in P}\frac{Q_c}{\prod_{j\in J,[D_j]=c}\alpha_j}}. 
\end{equation}

One observation is that $r$ of the factors $Q_c/(\prod_{j\in J,[D_j]=c}\alpha_j)$ are of degree $1$ as rational functions, and the others are of degree 0, since $\alpha_i, i\in \Sigma(1)-P-J$ form a basis of $W$. We denote the degree $1$ factors as $\tilde{\alpha}_i$, $i=1,...,r$, and the degree $0$ factors $l_j$. Rewrite $(\prod l_j) \tilde{\alpha}_1$ as $\tilde{\alpha}_1$ so that all $l_j$ are absorbed. After re-indexing, we can assume that $\tilde{\alpha}_i$ is a small deformation of $\alpha_i$. Then we claim that the integration result satisfies
\begin{equation}
|\int_{T_F(\epsilon)} \frac{\alpha_J}{\prod_{c\not\in P}Q_c} d\mu|=
|\int_{T_{\tilde{\alpha}}} \frac{1}{\prod_{i=1}^r \tilde{\alpha}_i} d\mu|,
\end{equation}
where $T_{\tilde{\alpha}} = \{u; |\tilde{\alpha}_j(u)/a_{jj}|=\epsilon_j\}$.
To see why this is true, we note that similarly to the above description of deforming $T_F$ to $T_\alpha$, we can deform $T_\alpha$ to $T_{\tilde{\alpha}}$, as long as $Q_c(\mathcal{E})$ is a sufficiently small deformation of $Q_c(T) = \prod_{i;[D_i]=c}\alpha_i$. Since $T_{\tilde{\alpha}}$ only differs by possibly an orientation from the permutation of $\alpha_i$ for different $F$'s, we conclude that the result is independent of the flag $F$ as long as 
 $\{\alpha_i, i\in \Sigma(1)-P-J\}$ forms a basis compatible with $F$. $\qed$

Now we are ready to prove Theorem \ref{$(0, 2)$classical}.
\subsection{Proof of Theorem \ref{$(0, 2)$classical}} 

Since Szenes-Vergne \cite{szenes2004toric} have proved the corresponding result in the $(2, 2)$ case, the map 
\begin{equation}
\sigma \mapsto \int_{T_F(\epsilon)} \frac{\sigma}{\prod_{c\in[\Sigma(1)]} Q_c}d\mu
\end{equation}
is not identically zero for small deformations. So it suffices to prove that for $\sigma\in SR(V,\mathcal{E})$ with $\deg \sigma=n$,
\begin{equation}\label{primitive0}
\int_{h(\mathcal{E})} \frac{\sigma}{\prod_{c\in[\Sigma(1)]} Q_c}d\mu =\sum_{F\in\mathcal{FL}^+(\xi)}\nu(F)\int_{T_F(\epsilon)} \frac{\sigma}{\prod_{c\in[\Sigma(1)]} Q_c}d\mu= 0.
\end{equation}

Since $SR(V,\mathcal{E})$ is generated by $$\{\prod_{c\in P}Q_c|P\subset \Sigma(1)  {\rm\ is \  a \ primitive \  collection }\}$$, it suffices to prove that the above equality (\ref{primitive0}) is true for those $\sigma$ of the form $\sigma=(\prod_{c\in P} Q_c)\cdot \sigma_J$, where $P$ is a primitive collection, $J\subset\Sigma(1)$, and $|J| =n -|P|$. 
Now we compute this by deforming the corresponding $(2, 2)$ result. 

Note that $\prod_{c\not\in P}Q_c$ is a small deformation of $\prod_{i\not\in P}\alpha_i$, we can write $\prod_{c\not\in P}Q_c = \prod_{i\not\in P}\alpha_i - \delta \tilde{\alpha}$, for some small $\delta\in\mathbb{C}$ and $\tilde{\alpha}\in \s^* W$. So we have 
\begin{equation}\label{power-series}
\begin{array}{ll}
\int_{h(\mathcal{E})} \frac{\sigma}{\prod_{c\in[\Sigma(1)]} Q_c}d\mu 
& = \int_{h(\mathcal{E})} \frac{\alpha_J}{\prod_{c\not\in P} Q_c}d\mu\\
& = \int_{h(\mathcal{E})} \frac{\alpha_J}{\prod_{i\not\in P}(\alpha_i - \delta \tilde{\alpha})}d\mu\\
& = \int_{h(\mathcal{E})} \frac{\alpha_J}{\prod_{i\not\in P}\alpha_i} (\sum_{n=0}^\infty (\frac{\delta\tilde{\alpha}}{\prod_{i\not\in P}\alpha_i})^n)d\mu\\
& = \sum_{m=0}^\infty \int_{h(\mathcal{E})} \frac{ (\delta\tilde{\alpha})^m\alpha_J}{(\prod_{i\not\in P}\alpha_i)^{m+1}}d\mu.
\end{array}
\end{equation}

Claim: for each monomial $\prod_{i\in K} \alpha_i$ such that $\frac{ \prod_{i\in K}\alpha_i}{(\prod_{i\not\in P}\alpha_i)^{m+1}}$ has degree $-r$,  
\begin{equation}
 \int_{h(\mathcal{E})}  \frac{ \prod_{i\in K}\alpha_i}{(\prod_{i\not\in P}\alpha_i)^{m+1}}d\mu =0.
 \end{equation} 

Proof of the Claim: by Lemma \ref{degenerate}, if the factors of the denominator do not generate $W$, the integration is 0. Hence it reduces to the case when the factors of the denominator generate $W$. 

For $k\in K$, write $\alpha_k = \sum a_{k_l}\alpha_l$, where $l$ runs through those indices appearing in the denominator. This reduces the integrand to
\begin{equation}
  \frac{ \prod_{i\in K}\alpha_i}{(\prod_{i\not\in P}\alpha_i)^{m+1}}d\mu = \sum_l \frac{a_{kl}\cdot\prod_{i\in K, i\neq k}\alpha_i}{ (\prod_{i\not\in P}\alpha_i)^{m+1} /\alpha_l} d\mu.
\end{equation}
Observe that as long as the remaining denominator of a summand generates $W$, we can repeat this process of expressing the numerator terms into linear combinations of the denominator terms and then canceling out a term. This process terminates after finitely many steps, and the final expression is a summation of terms of two types:

Type (i): terms with non-generating denominators. 

Type (ii): terms where factors of the denominator generate W, while the numerator is a constant (degree 0).

Type (i) terms integrate to 0 by Lemma \ref{degenerate}. So to prove the claim, it suffices to show that each Type (ii) term integrates to 0. Namely
\begin{equation}\label{-r}
\int_{h(\mathcal{E})} \frac{1}{\prod_{j=1}^r \alpha_{i_j}}d\mu=0.
\end{equation}

 Type (ii) terms have denominators of degree $r$, since the cancellation process preserves the degree of the fraction. Note that for the factors of the degree $r$ denominator to generate $W$ which is $r$ dimensional, these $r$ factors have to be distinct. So, being factors of $(\prod_{i\not\in P}\alpha_i)^{m+1}$, they are actually factors of $\prod_{i\not\in P}\alpha_i$.

Hence we have 
\begin{equation}
\begin{array}{ll}
&\int_{h(\mathcal{E})} \frac{1}{\prod_{j=1}^r \alpha_{i_j}}d\mu\\
=&\int_{h(\mathcal{E})} \frac{\prod_{i\in L}\alpha_i}{\prod_{i\in \Sigma(1)-P} \alpha_i}\\
=&\int_{h(\mathcal{E})} \frac{\prod_{i\in L\cup P}\alpha_i}{\prod_{i\in \Sigma(1)} \alpha_i}\\
=&0.
\end{array}
\end{equation}
The last equality comes from Theorem \ref{$(2, 2)$classical} and the fact that $\prod_{i\in L\cup P}\alpha_i$ is in $ SR(V)$. This proves Equation (\ref{-r}). Hence the claim is proved.

The theorem then follows from the claim and Equation (\ref{power-series}).

\section{The first integral formula for quantum correlators} \label{quantum-integral}
\subsection{The quantum correlators}\label{quantum-correlator}

Let $d_c^{\beta_j}$ be the intersection number of $\beta_j$ with any divisor in the divisor class $c$. 

The \textit{quantum correlator} of $\sigma_i\in W, i\in I$ is defined to be a summation over the GLSM moduli spaces $\mathcal{M}_\beta$ indexed by effective curves $\beta$:
\begin{equation}
\langle \sigma_I \rangle^{quantum} = \sum_\beta \langle \sigma_I F_\beta \rangle_\beta q^\beta,
\end{equation}
where $\beta$ runs over the lattice points in the Mori cone (generated by effective curve classes) of the toric variety $V$, and $F_\beta$ is the four-Fermi term introduced in \cite{donagi2014mathematical}:
\begin{equation}
F_\beta = \prod_{c\in [\Sigma(1)]} Q_c^{h^1(d_c^\beta)}.
\end{equation}
Also, when $\beta^{'}$ dominates $\beta$, the correlators over different moduli spaces are related by the ``exchange rate" $R_{\beta^{'}\beta}$:
\begin{equation}\label{exchangeRate}
\begin{array}{l}
\langle \sigma_I F_\beta \rangle_\beta = 
\langle \sigma_I F_\beta R_{\beta^{'} \beta} \rangle_{\beta^{'}}, \\
R_{\beta^{'}\beta} = \prod_c Q_c^{h^0(d^{\beta^{'}}_c)-h^0(d^{\beta}_c)}.
\end{array}
\end{equation}
\textbf{Remark:}
This holds even when $\beta$ is not effective.

We now define the complex value $q^\beta(z)$: For each $z=\sum_{i=1}^{n+r} z_i \omega_i\in g$ ($\omega_i$ corresponds to the ray $i\in\Sigma(1)$), and each $\beta\in H_2(V,\mathbb{Z})$, define $q^\beta(z)=\prod_{i=1}^{n+r}z_i^{\langle \alpha_i,\beta\rangle}$, where $\alpha_i \in H^2(V,\mathbb{Z})$ corresponds to the ray $i\in \Sigma(1)$. Sometimes we will write $q^\beta$ for $q^\beta(z)$.

\subsection{The Mori cone of a smooth projective variety}

The Mori cone of a variety is the closure of the cone effective curves. For projective toric variety $V$, the Mori cone $\overline{NE}(V)$ is a strongly convex rational polyhedral cone of full dimension in $N_1(V)$, the real vector space of proper 1-cycles modulo numerical equivalence. (\cite{cox2011toric}, pp 292 - 295.) Theorem 6.4.11 of \cite{cox2011toric} gives a concrete description of the Mori cone for any projective simplicial toric variety, representing it by specifying a generator for each primitive collection of $V$.\footnote{The simplical condition is not necessary, as commented in \cite{cox2011toric}, p 307.} As a corollary, we have:
\begin{prop}
The number of primitive collections of a simplicial projective toric variety $V$ is no less than its Picard number.
\end{prop}
Proof: By the above mentioned theorem, the number of primitive collections is the same as the number of cone generators of the Mori cone. Since the Mori cone has dimension $r =$ the Picard number, the proposition is proved. $\qed$
\subsection{The integral formula}\label{integral-formula-section}
We can choose $r$ generators $\beta_1,...,\beta_r$ such that they generates all the lattice points in a possibly bigger cone containing the Mori cone. We also require that $\langle \beta_j, \sum_{i=1}^{n+r} D_i\rangle \geq 0$. Since outside the Mori cone the moduli space is simply empty and (\ref{exchangeRate}) still holds in this case, we can write the summation over the Mori cone as summation over this (possibly bigger) cone:
\begin{equation}
\begin{array}{ll}
\langle \sigma_I \rangle^{quantum} 
& = \displaystyle \varinjlim_{B}\displaystyle\sum_{\beta {\rm \ dominated}{\rm \ by\ } B} 
\langle \sigma_I \cdot F_\beta \cdot R_{B\beta} q^{\beta} \rangle_B\\
& = \displaystyle \varinjlim_{B} \langle \sigma_I\displaystyle\sum_\beta \prod_c Q_c^{h^0(d^B_c) - h^0(d^\beta_c) + h^1(d^\beta_c) } q^\beta \rangle_B \\
& = \displaystyle \varinjlim_{B} \langle\sigma_I \prod_c Q_c^{h^0(d^B_c)-1} 
\displaystyle\sum_{\beta}  \prod_c Q_c^{-d_c^\beta} q^\beta \rangle_B \\
& = \displaystyle\lim_{N\to \infty} \langle\sigma_I \prod_c Q_c^{h^0(d^B_c)-1} 
\displaystyle \prod_{j=1}^r (\sum_{a_j=0}^{N} u_j^{a_j} )
\rangle_B\\
& = \displaystyle\lim_{N\to \infty} \langle\sigma_I \prod_c Q_c^{h^0(d^B_c)-1} 
\displaystyle \prod_{j=1}^r \frac{1 - u_j^{N+1}}{1-u_j}
\rangle_B,
\end{array}
\end{equation}

where  $u_j = \prod_c Q_c^{-d_c^{\beta_j}}q^{
\beta_j}$, and we will later write $q_j$ for $q^{\beta_j}$.

Now we have:
\begin{equation}
\begin{array}{ll}
&\langle \sigma_I \rangle^{quantum} \\
 = & \displaystyle\lim_{N\to \infty}  \frac{1}{(2\pi i)^r}\int_{h(\mathcal{E})} \frac{1}{\prod_c Q_c^{h^0(d^B_c)} }\cdot\left( \sigma_I \prod_c Q_c^{h^0(d^B_c)-1} 
\displaystyle \prod_{j=1}^r \frac{1 - u_j^{N+1}}{1-u_j}\right) d\mu\\
=&\displaystyle\lim_{N\to \infty}  \frac{1}{(2\pi i)^r}\int_{h(\mathcal{E})} \frac{1}{\prod_{c\in[\Sigma(1)]} Q_c} \cdot\left( \sigma_I 
\displaystyle \prod_{j=1}^r \frac{1 - u_j^{N+1}}{1-u_j}\right) d\mu\ .
\end{array}
\end{equation}
Let $Z(\mathbf{\epsilon})$ be a representative of $h(\mathcal{E})$ and $q_j$ be sufficiently small. We have $|u_j| < 1$ on $h(\mathcal{E})$. We also write $v_j = u_j^{-1} =\prod_c Q_c^{d_c^{\beta_j}}q_j^{-1}$, $\tilde{v}_j=\prod_c Q_c^{d_c^{\beta_j}}$. So we have $|v_j|>1$ on $h(\mathcal{E})$. Hence
\begin{equation}
\begin{array}{ll}
\langle \sigma_I \rangle^{quantum}  &=\displaystyle\lim_{N\to \infty}\frac{1}{(2\pi i)^r}
 \int_{Z(\mathbf{\epsilon})} \frac{\sigma_I }{\prod_c Q_c } \prod_{j=1}^r \frac{1 - u_j^{N+1}}{1-u_j}
d\mu \\
&=\displaystyle\lim_{N\to \infty}\frac{1}{(2\pi i)^r}
 \int_{Z(\mathbf{\epsilon})} \frac{\sigma_I }{\prod_{c\in[\Sigma(1)]} Q_c} \prod_{j=1}^r \frac{1 - u_j^{N+1}}{1-u_j}
d\mu \\
& =\displaystyle \frac{1}{(2\pi i)^r}
 \int_{Z(\mathbf{\epsilon})} \frac{\sigma_I }{\prod_c Q_c } \prod_{j=1}^r \frac{1 }{1-u_j}
d\mu\\
& = \displaystyle  \frac{1}{(2\pi i)^r}
 \int_{Z(\mathbf{\epsilon})} \frac{\sigma_I \prod_j v_j }{\prod_{c\in[\Sigma(1)]} Q_c} \prod_{j=1}^r \frac{1 }{v_j-1}
d\mu
\end{array} 
\end{equation}

Thus we have proved the following result:
\begin{thm}\label{integral-formula}
Let $\mathcal{E}$ be a holomorphic vector bundle defined by the deformed toric Euler sequence (\ref{euler1}) with small deformations.
Let $Z(\mathbf{\epsilon})$ be a cycle representing $h(\mathcal{E})$ and $z\in (\mathbb{C}^*)^n$. Let $q_j = z^{\beta_j}$, $j=1,...,r$. For a fixed basis $\beta_1,...,\beta_r$ and $z\in (\mathbb{C}^*)^n$ such that $|q_j| < \min_{u\in Z(\mathbf{\epsilon})}|\tilde{v}_j(u)|$ holds, we have 
\begin{equation}\label{LHS}
\langle \sigma_I\rangle^{quantum} = \displaystyle  \frac{1}{(2\pi i)^r}\displaystyle\int_{Z(\mathbf{\epsilon})} \frac{\sigma_I}{\prod_{c\in [\Sigma(1)]}Q_c} \frac{\prod_{j=1}^{r}\tilde{v}_j}{\prod_{j=1}^r (\tilde{v}_j -q_j)}d\mu.
\end{equation}
\end{thm}
\section{Number of solutions to the QSC Relations}\label{QSC-system}

The QSC relations $\tilde{v}_j =q_j$ form a Laurent polynomial system $\mathcal{F}_0$. In order to prove the main result of the paper, we need to study the solutions to the system $\mathcal{F}_0$.

\subsection{The toric case}

In \cite{donagi2015global-aspects} the QSC for toric deformations \footnote{meaning that the deformed bundle is toric equivariant.} of the tangent bundle are studied. We briefly recall the relevant results here.

In this case the QSC can be written as: 

\begin{equation}
\prod_\rho (\sum_{b=1}^r E^b_\rho \sigma_b)^{d^a_\rho} = q_a,
\end{equation}

where $q_a\in \mathbb{C}^*$ and $(E^a_\rho)$ is of maximal rank.

Also note that for the tangent bundle, $E^a_\rho=d^a_\rho$.

Let $z_\rho = \sum_{b=1}^r E^b_\rho \sigma_b$. Then the $z_\rho$'s satisfy $ \hat{E} \cdot z = 0$,
where $\hat{E}$ is a matrix corresponding to the Gale dual of $E$, i.e. the $n\times (n+r)$ matrix whose rows record the linear relations of the columns of $E$. 

In this way the system of QSC for toric deformations of the tangent bundles can be written as 
\begin{equation}\label{QSCRfull}
\begin{array}{l}
\prod_\rho z_\rho^{d^a_\rho} = q_a,\\
\hat{E} \cdot z = 0.
\end{array}
\end{equation}

Write the $i$-th equation as $F_i = 0$ and denote its Newton polytope by $S_i$. Let $\mathcal{F} =(F_1,...,F_N)$ and $S = (S_1,...,S_N)$.

By \cite{cox2011toric} Proposition 6.4.1, we have the following short exact sequences:
\begin{equation}\label{ses_N_1}
\begin{array}{l}
0 \to M_\mathbb{R} \xrightarrow{\alpha} \mathbb{R}^{\Sigma(1)} \xrightarrow{\beta} H^2(X) \to 0,\\
0 \to N_1(X)_\mathbb{R} \xrightarrow{\beta^*} \mathbb{R}^{\Sigma(1)} \xrightarrow{\alpha^*} N_{\mathbb{R}} \to 0,
\end{array}
\end{equation}
where

\begin{tabular}{ll}
$\alpha^*(e_\rho) =u_\rho$, & $u_\rho$ the primitive vector in the ray $\rho$.\\
$\beta^*([C]) = (D_\rho\cdot C)_{\rho\in\Sigma(1)}$, & $C\subset X$ an irreducible complete curve.
\end{tabular}

Thus it is clear that $S_a, a=1,...,r$ spans $\beta^*(N_1(X))$, and its orthogonal complement $N_1(X)^\perp$ is isomorphic to $N_{\mathbb{R}}$. Also, $S_{r+1}=...=S_N = $ Conv$(e_1,...,e_N)$.

By Bernstein \cite{bernstein1975number}, for generic $(E,q)$, the number of solutions to this Laurent polynomial system (QSC system for short) in $(\mathbb{C}^*)^n$ is the mixed volume of $S$. This number is $\chi(X)$ by \cite{donagi2015global-aspects} Theorem.... If a Laurent polynomial system with coefficients $(E,q)$ has the desired number of solutions, we will call it \textit{Bernstein generic}.

\subsection{The $(2, 2)$ case}

Properness in the $(2, 2)$ case means that for $q\in R(\xi,\eta)$, under some technical conditions, the Laurent polynomial system has finitely many solutions, and the number is at most $\chi(X)$.

\begin{thm}\label{TX_is_generic}
Let $X$ be a smooth projective nef-Fano toric variety. For the tangent bundle $T_X$ and generic $q$, $(T_X,q)$ is Bernstein generic in QSC for toric deformations.
\end{thm}

Proof: By Bernstein, one needs to check that for every $w\neq 0$, there is no solution to $F^w$ for generic $q\in(\mathbb{C}^*)^N$. 

Note that for $ i\leq r$, $S_i$ is one-dimensional. Hence for $F^w$ to have solutions with $q\in(\mathbb{C}^*)^N, a=1,...,r$, it is necessary to have $w\perp S_i, i\leq r$, i.e. for any $v\in N_1(X)$, $\langle \beta w,v \rangle = \langle w, \beta^* v\rangle =0$. Hence there exist $w_0\in M_\mathbb{R} $ such that $w = \alpha w_0$. Hence $\langle w,e_\rho \rangle =\langle \alpha w_0, e_\rho \rangle =\langle w_0, \alpha^* e_\rho\rangle = \langle w_0, u_\rho\rangle$. Hence $F^w$ is completely determined by the corresponding $w_0$ and the bilinear pairing of $N_\mathbb{R}$ and $M_\mathbb{R}$.

In Fano case, this is easy. For $w_0\neq0$, $\alpha^* (S^w_N)$ is always of the form Conv$(u_\rho, \rho\in \sigma(1))$, for some $\sigma\in\Sigma$. Hence the system $F^w_j = 0, j=r+1,...,N$ has no nonzero solution because the non-degeneracy condition on the corresponding columns of the matrix $\hat{E}$ for $T_X$. This means that in Fano case, $(T_X,q)$ is Bernstein generic for any $q$.

When $X$ is nef-Fano but not Fano, by the previous argument, in order for $(T_X, q)$ to be non Bernstein generic, it is necessary that there exist $w\neq 0$ and $m>0$, such that 
$\alpha^*(S^w_N) =\langle u_1,...,u_{n^\prime+m}\rangle$, where $n^\prime =\dim \alpha^*(S^w_N) + 1 \leq n$. In this case we can write the solutions to $F^w_{r+1} = ... = F^w_{N} = 0$ as 
$z_j = z_{n^\prime +1}\cdot L_j(\frac{z_{n^\prime+2}}{z_{n^\prime +1}},...,\frac{z_{n^\prime+m}}{z_{n^\prime +1}})$, where $L_j$ is a linear function, $j=1,...,n^\prime$, up to a permutation of the indices in $\{1,...,N\}$.

Now consider the $n^\prime + m$ rays corresponding to $u_1,...,u_{n^\prime + m}$ in a $n^\prime$ dimensional subspace of $N_\mathbb{R}$. There are $m$ independent linear relations between them, giving rise to $q_a = \prod z_\rho^{d^a_\rho}$, $a =1,...,m$, with $\sum_{\rho =1}^{n^\prime + m} d^a_\rho = 0$, and $d^a_\rho =0$, $\rho > n^\prime + m$. 

This implies that $q_1,...,q_m$ depends only on $m-1$ free variables $\frac{z_{n^\prime+2}}{z_{n^\prime +1}},...,\frac{z_{n^\prime+m}}{z_{n^\prime +1}}$. Hence non Bernstein generic $(T_X,q)$'s are not generic. This finishes the proof.\qed

\subsection{The general case}

We say a ray $\rho\in \Sigma$ (and the corresponding toric divisor $D_\rho$) is \textit{rigid}, if there is no other $\rho^\prime\in\Sigma$ such that $D_\rho \sim D_{\rho^\prime}$. Otherwise we say $\rho$ and $D_\rho$ are non-rigid. In other words, let $n_c$ be the number of prime toric divisors in the divisor class $c$. Then $D_\rho$ is rigid iff $n_{[D_\rho]}=1$. We will denote $n_{D_\rho}$ by $n_\rho$ for short.

For linear deformations of the tangent bundle, the QSC relations can be written as 
\begin{equation}\label{N_plus_gamma_equations}
\begin{array}{rl}
\displaystyle\prod_{\rho\in\{\rho | n_\rho=1\}}z_\rho^{d_\rho^a}\cdot \prod_{c\in\{c | n_c >1\}} y_c^{d_c^a}  &= q_a,\\
\hat{E} \cdot z & =0,\\
y_c - Q_c(z_1,...,z_N) & =0,
\end{array}
\end{equation} 
which is a deformed version of (\ref{QSCRfull}). Note the abuse of notation; here we have made a choice of $\sigma_a=\sigma_a(z_1,...,z_N)$ subject to $z_\rho = \sum E^b_\rho \sigma_b$, which makes $Q_c(\sigma_1,...,\sigma_r)=Q_c(z_1,...,z_N) \in \mathbb{C}[z_1,...,z_N]$. Note that the construction of the QSC implies that each $Q_c(z_1,...,z_N)$ is a homogeneous polynomial of degree $n_c$. 

We can form two Laurent polynomial systems with variables  $$ (z_1, ..., z_N, y_1, ..., y_\gamma) $$ based on (\ref{N_plus_gamma_equations}) and different choices of coefficients of $Q_c$. For the first one, we let $Q_c = \prod_{\rho\in\{\rho|[D_\rho]=c\}} z_\rho$, which corresponds to the QSC for toric deformations. For the second one, let $Q_c =\sum_{|I|=n_c} a_I z^I$, which is the most general form of homogeneous polynomial of degree $n_c$. This covers all possible linear deformations.
Denote the two systems as $\mathcal{F} = (F_1,...,F_{N+\gamma})$ and $\tilde{\mathcal{F}} = (\tilde{F}_1,...,\tilde{F}_{N+\gamma})$, and the corresponding Newton polytopes as $S=(S_1,...,S_{N+\gamma})$ and $\tilde{S} = (\tilde{S}_1,...,\tilde{S}_{N+\gamma})$, where $\gamma = |\{c|n_c>1\}|$. Let $MV(\mc{S})$ be the mixed volume of $\mc{S}$. Note that the newly-defined $\mathcal{F}$ is equivalent to the system (\ref{QSCRfull}) with $MV(\mc{S})=\chi(X)$, which is the reason why we did not introduce a new notation.

\begin{thm}\label{general_deformation}
For any smooth projective nef-Fano toric variety $X$ and generic $q$, a generic linear deformation of the tangent bundle preserves the number of solutions to the QSCR, which is $\chi(X)$.
\end{thm}

This is an important technical ingredient to the proof of the properness lemma in next section (Lemma \ref{properness}).

We need to use the concept of \textit{essential subsets} to prove this theorem. 

\begin{defn}(\cite{rojas1994convex})
Let $(C_1,...,C_n)$ be an $n-$tuple of polytopes or an $n$-tuple of finite point sets in $\mathbb{R}^n$. A non-empty subset $I\subset \{1,...,n\}$ is said to be an \textit{essential subset} for $(C_1,...,C_n)$ if $\dim \sum_{i \in I} C_i = |I|-1$ and $\dim \sum_{i\in J} C_i \geq |J|-1$ for all proper nonempty $J\subset I$. 
\end{defn}

Fact: $(C_1,...,C_n)$ has an essential subset $\Leftrightarrow$ $MV(C_1,...,C_n) = 0$. 

\begin{lem}(\cite{rojas1994convex}, Corollary 9)\label{MVEquality}
Suppose $\mathcal{S} = (S_1,...,S_n)$ and $\tilde{\mathcal{S}} = (\tilde{S}_1,...,\tilde{S}_n)$ are $n$-tuples of polytopes in $\mathbb{R}^n$ with rational vertices, such that $MV(\tilde{\mathcal{S}}) > 0$ and $\mathcal{S}\subset \tilde{\mathcal{S}}$. Then $MV(\mathcal{S}) =MV(\tilde{\mathcal{S}})$ $\Leftrightarrow$ for all $w\in \mathcal{S}^{n-1}$, $\tilde{\mathcal{S}}^w$ has an essential subset $I^w$ such that $\tilde{S}^w_i \cap S_i \neq \phi$ for all $i\in I^w$. 
\end{lem}
Proof of Theorem \ref{general_deformation}: In the language of Bernstein, it suffices to prove that the mixed volumes of $\tilde{S}$ and $S$, $MV(\tilde{\mc{S}})$ and $MV(\mc{S})$ are equal. 

Since $S_i\subset \tilde{S}_i$, we have $MV(\tilde{\mc{S}}) \geq MV(\mc{S})$.

The strategy to proceed is to use the separation lemma.

We have $\tilde{S}_i \subset \mathbb{R}^{N+\gamma}$. We will use the short exact sequence 
\begin{equation}
0 \to N_1(X)_\mathbb{R} \xrightarrow{\tilde{\beta}^*} \mathbb{R}^{N+\gamma} \xrightarrow{\tilde{\alpha}^*} \widetilde{N_\mathbb{R}}\to 0.
\end{equation}

Let $A^\prime= \tilde{\alpha}^*(A)$, for any $A\subset \mathbb{R}^{N+\gamma}$. Note that for $i=1,...,r$, $S_i =\tilde{S}_i \subset \tilde{\beta}^*(N_1(X)_\mathbb{R})$. So it suffices to prove that $$MV(\tilde{S}^\prime_{r+1},...,\tilde{S}^\prime_{N+\gamma})= MV(S^\prime_{r+1},...,S^\prime_{N+\gamma}).$$

Since $S_i^\prime\subset \tilde{S}^\prime_i$, it suffices to show that for any $w\neq 0$, there is an essential subset $J$, such that for any $j\in J$, $(\tilde{S}^\prime_j)^w \cap S^\prime_j \neq \phi$ by Lemma \ref{MVEquality}.

Notice that $w$ is in the normal fan of the polytope $\sum_{i=r+1}^{N+\gamma} \tilde{S}^\prime_i$. If $w$ and $\tilde{w}$ are in the relative interior of the same cone, $(\tilde{S}^\prime_j)^w =(\tilde{S}^\prime_j)^{\tilde{w}}$. Hence there are only finitely many $w$'s to check.

Since $Im\ \tilde{\beta}^* \subset (\mathbb{R}^{k+\gamma})^\perp$, we have $\widetilde{N_\mathbb{R}}\cong \overline{N_\mathbb{R}} \oplus \mathbb{R}^{k+\gamma}$. We introduce some notations and express $S^\prime_{j}$ and $\tilde{S}^\prime_{j}$ using them. 

Let $w= (w^\prime, w^{\prime\prime})$, where $w^\prime \in \overline{N_\mathbb{R}}^*$ and $w^{\prime\prime} \in (\mathbb{R}^{k+\gamma})^*$. We write $w^{\prime\prime} = (w_1,...,w_{k+\gamma})$, 
with $w_j$ corresponding to the non-rigid ray $j$. Let $w^{\prime\prime}_{\min} = \min\{w_j | j \ \text{non-rigid}\}$.
Note that under $\tilde{\alpha}^*$, we have 
$$e_j\mapsto \left\{\begin{array}{ll}
{u}_j,& j<N, j \ \text{rigid},\\
e_j,& j<N, j \ \text{non-rigid},\\
v_\delta=\sum_{l\in K_\delta} u_l, & j = N+\delta.
\end{array}\right.$$
Let $\bar{v}_\delta =\frac{1}{k_\delta} v_\delta$. Let $M = \{j|w_j = w^{\prime\prime}_{\min}\}$ and $m = |M|$.
We then have
$\tilde{S}^\prime_{r+j} = S^\prime_{r+j} = Conv({u}_\rho, e_\tau|\rho,\tau\in\Sigma(1), \rho \ \text{rigid}, \tau \ \text{non-rigid})$, $j = 1, 2, ..., n$. 
$S^\prime_{N+\delta} = Conv(v_\delta,\sum_{l\in K_\delta} e_l), \tilde{S}^\prime_{N+\delta} = Conv(v_\delta,\sum_{l\in I} e_l)_{|I|=|K_\delta|}$.

Now we proceed by considering all the different cases of $w$:

(1) $\displaystyle\min_{j\ \rm{rigid}} \langle w, {u}_j \rangle  < w^{\prime\prime}_{\min}$: then $ \tilde{S}^{\prime w}_{r+1} = ... =  \tilde{S}^{\prime w}_N  \subset \overline{N_\mathbb{R}}$. Note that we always have $\tilde{S}^{\prime w}_j \cap S^{\prime}_j \neq \phi, j=r+1,...,N$. And we can always pick an essential subset among them for dimension reason: let $m =\dim \tilde{S}^{\prime w}_N \leq n - k$, with $k = 0 $ being the trivial case where there are no non-rigid rays. For $k>0$, we have $m\leq n-1$, and any $J\subset \{r+1,...,N\}$ with $|J| = m+1$ is an essential subset as required.

(2) $\min \langle w, \bar{u}_j \rangle  \geq w^{\prime\prime}_{\min}$. There are two sub-cases:

(2.1) $w^\prime = 0$, hence $\langle w, {u}_j \rangle =0$ for any $j$. There are two sub-cases:

(2.1.1) $w^{\prime\prime}_{\min} = 0$. In this case $\langle w, v_\delta \rangle =0= k_\delta  w^{\prime\prime}_{\min}\leq $ other elements of $\tilde{S}_{N+\delta}^\prime$, for any $\delta$. Hence $v_\delta \in\tilde{S}^{\prime w}_{N+\delta} \cap S^{\prime w}_{N+\delta}$, for any $\delta$. Hence any essential subset satisfies the requirement.

(2.1.2) $w^{\prime\prime}_{\min} < 0$: In this case $\tilde{S}^{\prime w}_j = \tilde{S}^{\prime w}_N, j = r+1,...,N$ (NOT $N+\gamma$).  If $m \leq k$ (hence $m\leq n$), any $J\subset \{r+1,...,N\}$ with $|J| = m$ is an essential subset since $\dim \tilde{S}^{\prime w}_N = m -1$. If $m > k$, then there are at least $(m-k)$ $\delta$'s such that $K_\delta\subset M$ (pigeon hole principle). Denote them as  $K_1,...,K_{m-k}$, then $\{r+1,...,r+k\} \cup \{N+1,...,N+m-k \}$ is an essential subset as required.

(2.2) $w^\prime \neq 0$, hence won't be able to pick up all $\bar{u}_j$'s in $\tilde{S}^{\prime w}_{j\leq N}$. There are two sub-cases:

(2.2.1) There exists $\delta$ such that $\langle w, \bar{v}_\delta\rangle  <  w^{\prime\prime}_{\min}$, then $\{N+\delta\}$ is an essential subset with required property.

(2.2.2) $\langle w, \bar{v}_\delta\rangle  \geq  w^{\prime\prime}_{\min}$. In this case for $j=r+1,...,N$, $e_i\in \tilde{S}^{\prime w}_j$ iff $ i\in M$. Hence $\dim \tilde{S}^{\prime w}_j \leq n-k-1 + m$, $j=r+1,...,N$. Moreover, $\dim \sum_{j\in J}\tilde{S}^{\prime w}_j \leq n - k - 1 + m$, if $J\subset \{r+1,...,N\}$. If $m\leq k$, $\dim \sum_{j\in J}\tilde{S}^{\prime w}_j \leq n-1$, we can easily pick an essential subset $J\subset \{r+1,...,N\}$ with required property. If $m > k$, there are at least $(m-k)$ $\delta$'s such that $K_\delta\subset M$. Denote them as  $K_1,...,K_{m-k}$, then $J_2=\{r+1,...,N+m-k \}$ contains an essential subset with required property. The reason is that we still have $\dim \sum_{j\in J}\tilde{S}^{\prime w}_j \leq n-1$. To see this, note that $\langle w, u\rangle \equiv w^{\prime\prime}_{\min}$ and $w^\prime \neq 0$. So we know that $\tilde{S}^{\prime w}$ cannot contain the full $\overline{N}_{\mathbb{R}}^*$. Hence we still have  $\dim \sum_{j\in J_2}\tilde{S}^{\prime w}_j \leq n-1$.

In conclusion, we have shown that the essential set $J$ with required property exists in every case, which completes the proof.
\qed

\section{Quantum correlator: summation formula}\label{quantum-sum}

\subsection{The Main Result}\label{section_main_result}
In McOrist-Melnikov \cite{mcorist2008half}, there is a summation formula for quantum correlators, for $\mathcal{E}$ defined by the deformed toric Euler sequence (\ref{euler1}) with a linear deformation, as stated in Conjecture \ref{toric_case}. The authors derive it from physics argument using Coulomb branch techniques. We now prove it using quantum sheaf cohomology. For technical reasons we assume $V$ is \textit{nef-Fano}, namely, the anti-canonical class of $V$ is nef.

\begin{mr}
Let $V$ be a nef-Fano smooth projective toric variety.
Let $\mathcal{E}$ be a holomorphic vector bundle defined by the deformed toric Euler sequence (\ref{euler1}) with small deformations. Then 
\begin{equation}\label{RHS}
\langle \sigma_{i_1},...,\sigma_{i_s}\rangle^{quantum} =\displaystyle \sum_{\{u\in W^\vee|\tilde{v}_j(u)=q_j\}} \frac{\sigma_I}{\prod_{c\in [\Sigma(1)]} Q_c} \frac{\prod_{j=1}^{r}\tilde{v}_j}{\det_{j,k} (\tilde{v}_{j,k})}
\end{equation}
holds for $z$ in a non-empty open region in $(\mathbb{C}^*)^n$. Here terms on the right hand side, being elements in $\Sym^*W$, are viewed as functions on $W^\vee$.

\end{mr}


We explain the notations.  $f_{,k}$ stands for $\frac{\partial f(u)}{\partial u_k}$. Recall from Section \ref{integral-formula-section} that  $$\tilde{v}_j=\prod_{c\in [\Sigma(1)]} Q_c^{d_c^{\beta_j}},$$ where $z\in (\mathbb{C}^*)^n$, $q_j = z^{\beta_j}$, $j=1,...,r$. $\beta_1,...,\beta_r$ generates a cone containing the Mori cone. Let $z\in (\mathbb{C}^*)^n$ such that $|q_j| < \min_{u\in Z(\mathbf{\epsilon})}|\tilde{v}_j(u)|$ holds. 

\textbf{Remarks:}


(a) Let $h_q\in H_r(W^\vee-\{u\in W^\vee\ |\ \tilde{v}_j(u)=q_j\},\mathbb{Z})$ be the homology of the real $r$-dimensional cycle defined by $\{u\in W^\vee\ |\ |\tilde{v}_j(u) -q_j|< \delta\}$, for a $\delta$ that is small enough. Then the above formula (\ref{RHS}) can be written as 
\begin{equation}\label{int-hq}
\langle \sigma_I\rangle^{quantum} =
\displaystyle  \frac{1}{(2\pi i)^r} \displaystyle\int_{h_q} \frac{\sigma_I}{\prod_{c\in [\Sigma(1)]}Q_c} \frac{\prod_{j=1}^{r}\tilde{v}_j}{\prod_{j=1}^r (\tilde{v}_j -q_j)}d\mu.
\end{equation}

Equation (\ref{LHS}) and (\ref{int-hq}) has exactly the same integrand on their right-hand sides.


(b) When $\mathcal{E}$ is the tangent bundle, (\ref{int-hq}) is reduced to 
\begin{equation}
\langle \sigma_I\rangle^{quantum} =
\displaystyle  \frac{1}{(2\pi i)^r}\displaystyle\int_{h_q} \frac{\sigma_I}{\prod_{i\in\Sigma(1)}\alpha_i} \frac{\prod_{j=1}^{r}\tilde{v}_j}{\prod_{j=1}^r (\tilde{v}_j -q_j)}d\mu.
\end{equation}

This resembles the conclusion of the hypersurface case of the ``Toric Residue Mirror Conjecture" in the $(2, 2)$ theory \cite{batyrev2002toric}, which says\footnote{We adopt the formulation of  \cite{szenes2004toric}. See Proposition 4.7 there.} for anti-canonical hypersurface $X$ (i.e. the fundamental class is dual to $\kappa =\sum_{i\in \Sigma(1)}\alpha_i$) in a Fano simplicial toric variety $V$ of dimension $n$, the quantum correlator $\langle\sigma_{i_1}...\sigma_{i_{n-1}}\rangle^{quantum}$ for $\sigma_i\in H^1(X, T_X^*)$ coming from the restriction of $H^1(V,T_V^*)$ is
\begin{equation}\label{trmc}
\langle \sigma_{i_1}...\sigma_{i_{n-1}}\rangle^{quantum} = \displaystyle  \frac{1}{(2\pi i)^r}\displaystyle\int_{h_q} \frac{\sigma_I}{(1-\kappa)\prod_{i\in\Sigma(1)}\alpha_i} \frac{\prod_{j=1}^{r}\tilde{v}_j}{\prod_{j=1}^r (\tilde{v}_j -q_j)}d\mu.
\end{equation}

(c) The nef-Fano constraint is imposed since the proof depends on an alternative description of $Z(\epsilon)$ which only holds in the nef-Fano case. See 
Corollary \ref{z-represents-h_c} below. Results in previous sections do not require $V$ to be nef-Fano.

The main goal of this section is to prove the Main Result by proving (\ref{int-hq}). Our proof is inspired by Szenes and Vergne's proof of (\ref{trmc}).

We set up some notations in Section \ref{preparations}. Then we state the theorem and some lemmas. We prove the theorem in Section \ref{proof}.
\subsection{Some preparations}\label{preparations}

We fix a bundle $\mathcal{E}_1$ such that Theorem \ref{integral-formula} holds. Multiply the deformation parameters in the map defined $\mathcal{E}_1$ by $t$, we get a one parameter family $\mathcal{E}_t$. Then Theorem \ref{integral-formula} holds for $|t|\leq 1$ and in particular $\mathcal{E}_0 = T_V$.

Specifying the $t$-detpendence of the map $\tilde{v}_j$, we denote it by $\tilde{v}_j^{t}$.

We make some definitions generalizing those $(2, 2)$ case notations in \cite{szenes2004toric} to the $(0, 2)$ case: 

Define \begin{equation}
\hat{Z}^{t}(\xi)=\{u\in U; |\tilde{v}_j^{t}| = e^{-\langle \xi,\beta_j\rangle}\}.
\end{equation}
$\hat{Z}^{t}(\xi)$ can be viewed as the preimage of a torus $$T(\xi) =\{y\in (\mathbb{C}^*)^r; |y_j| = e^{-\langle \xi,\beta_j\rangle}\},$$ under the map $\tilde{v}^{t}=(\tilde{v}_1^{t},...,\tilde{v}_r^{t}): U(\mathcal{E}) \to (\mathbb{C}^*)^r$.

For $S\subset\{1,2,...n\}$ define
$$T_S(\xi,\eta)=\left\{ y\in (\mathbb{C}^{*})^r; |y_j|=
\left\{\begin{array}{l}
\exp(-\langle \xi,\beta_j\rangle ),{\rm\ if\ } j\in S,\\
\exp(-\langle \xi-\eta,\beta_j\rangle ),{\rm\ if\ } j\not\in S.
\end{array} \right.\right\}$$ 
and $T_\delta(q) =\{y\in (\mathbb{C}^*)^r; |y_j-q_j|=\delta, j=1,...,r\}.$ Let $Z_S^{t}(\xi,\eta)$ and $Z^{t}_\delta(q)$ be the pull-back of $T_S(\xi,\eta)$ and $T_\delta(q)$ respectively by $\tilde{v}^{t}$. Note that results about $\hat{Z}(\xi)$ apply to $\hat{Z}_S^{t}(\xi,\eta)$.

In order to keep the notations clean, we omit the label $t$ when $t=1$, 
and simply write $Z_S, Z_\delta, \tilde{v}$. Note that $Z_\delta(q)$ represents $h_q$.

Let $R(\xi,\eta)$ be the multi-dimensional annulus
\begin{equation*}
\{y = (y_1,...,y_r)\in (\mathbb{C}^*)^r; \langle \xi-\eta, \beta_j\rangle < -\log |y_j| < \langle\xi,\beta_j\rangle,j=1,...,r  \},
\end{equation*}
and let 
$W(\xi, \eta)$ be the pull-back of $R(\xi,\eta)$ by the map $q(z)$.

To simplify notation, define $$\Lambda = \displaystyle  \frac{1}{(2\pi i)^r}\displaystyle  \frac{\sigma_I}{\prod_{c\in [\Sigma(1)]}Q_c} \frac{\prod_{j=1}^{r}\tilde{v}_j}{\prod_{j=1}^r (\tilde{v}_j -q_j)}d\mu .$$ 

We also recall the definition of $\tau$-regularity from \cite{szenes2004toric}:

\begin{defn}(\cite{szenes2004toric})
$\mathfrak{U} = \{\alpha_i\in W; i\in \Sigma(1) \}$. Define 
\begin{equation*}
\Sigma\mathfrak{U} =\left\{ \sum_{i\in \eta} \alpha_i; \eta\subset \Sigma(1)\}\right\},
\end{equation*}
which is the collection of partial sums of elements of $\mathfrak{U}$. For each subset $\mathbf{\rho}\subset \Sigma\mathfrak{U}$ which generates $W$, we can write $\xi=\sum_{\gamma\in\mathbf{\rho}}a_\gamma^\mathbf{\rho}(\xi)\gamma$. Denote 
\begin{equation*}
\min(\Sigma\mathfrak{U},\xi) = \min\{|a_\gamma^\mathbf{\rho}(\xi)|; \mathbf{\rho}\subset\Sigma\mathfrak{U},\mathbf{\rho}{\rm \ basis \ of\ } W,\gamma\in\mathbf{\rho} \}.
\end{equation*}
We say $\xi\in W$ is $\tau$-regular for $\tau>0$ if $\min(\Sigma\mathfrak{U},\xi) >\tau$.
\end{defn}

The requirement for the technical assumption for $\xi$ being $\tau$-regular for sufficiently large $\tau$ will be seen in Proposition \ref{smooth-compact-cycle}.

Let $\xi$ be $\tau$-regular for $\tau$ sufficiently large. In the $(2, 2)$ case,  $\tilde{v}^{(0)}$ is regular over $(\tilde{v}^{(0)})^{-1}(R(\xi,\eta))$. Hence, as a consequence, we have

\begin{prop}\label{regularity of v}
 $\tilde{v}^{t}$ is regular over $(\tilde{v}^{(0)})^{-1}(R(\xi,\eta))$.
\end{prop}

\begin{lem}[The properness lemma]\label{properness}
Let $\xi$ be $\tau$-regular for $\tau$ sufficiently large.
Then the map $\tilde{v}^{t}$ is proper from $(\tilde{v}^{t})^{-1}(R(\xi,\eta))$ to $R(\xi,\eta)$ for all $\mathcal{E}$ sufficiently close to $T_X$ (with respect to the metric induced by the classical topology).
\end{lem}


Proof: the $(2, 2)$ case is proved in \cite{szenes2004toric}\footnote{It is Proposition 5.15, the map $\tilde{v}^{(0)}$ is just the map $p$ in \cite{szenes2004toric}.}. Theorem \ref{TX_is_generic} and Theorem \ref{general_deformation} shows that $T_X$ is Bernstein generic for $q\in R(\xi,\eta)$. Hence $\tilde{v}^{t}$ is proper from $(\tilde{v}^{t})^{-1}(R(\xi,\eta))$ to $R(\xi,\eta)$ for $\mathcal{E}$ sufficiently close to $T_X$.\qed
~\\


The technical assumption for $\xi$ being $\tau$-regular for sufficiently large $\tau$ is made to achieve the following result in the $(2, 2)$ case:
\begin{prop}(Theorem 6.2 of \cite{szenes2004toric})\label{smooth-compact-cycle}
If $\tau$ is sufficiently large, then for any $\tau$-regular $\xi\in\mathfrak{c}$, the set $\hat{Z}^{(0)}(\xi)$ is a smooth compact cycle whose homology class equals $h(\mathcal{E})\in H_r(U(\mathcal{E}_0),\mathbb{Z})=H_r(U(T_V),\mathbb{Z})$.
\end{prop}

\begin{cor}\label{z-represents-h_c}
The homology class of $\hat{Z}^{t}(\xi)$ is $h(\mathcal{E}_t)\in H_r(U(\mathcal{E}),\mathbb{Z})$.
\end{cor}

Proof: Lemma \ref{properness} shows the properness of $\tilde{v}^{t}$ on $(\tilde{v}^{t})^{-1}(R(\xi^\prime,\eta^\prime))$.

In order to make a homotopy argument, we first show that there is a $T > 0$ such that $\tilde{v}^{t_1}((\tilde{v}^{t_2})^{-1}(R(\xi,\eta))) \neq 0$, for any $t_1, t_2 \in \Delta_T$, where  $\Delta_T= \{t\in\mathbb{C};|t<T|\}$.

Observe that $\tilde{v}^{-1}(\overline{\Delta_T\times R(\xi,\eta)})$ is compact. Since $Q_i(t,u) \neq 0$ on $\tilde{v}^{-1}(\overline{\Delta_T\times R(\xi,\eta)})$, we can define $$\epsilon_i = \displaystyle\min( Q_i(t,u) | (t,u)\in\tilde{v}^{-1}(\overline{\Delta_T\times R(\xi,\eta)})$$ and have $\epsilon_i >0$.

Now if $\tilde{v}(t_2,u)\in R(\xi,\eta)$, then \begin{equation}
|Q_i(t_1,u) - Q_i(t_2,u)| = |t_1 - t_2| \cdot |\sum_I c_I z^I| < T^\prime M.
\end{equation}
Let $T^\prime = \min(T,\frac{\epsilon_i}{M} |i=1,2,...,r)$, then $Q_i(t_2, u ) \neq 0$. We rename $T^\prime$ to $T$ from now on. Then for $t_1,t_2 \in \Delta_T $, we have $\tilde{v}^{t_1}((\tilde{v}^{t_2})^{-1}(R(\xi,\eta))) \neq 0$.

Next, since $\tilde{v}^{t}$ is regular on this region, the cycles $\hat{Z}^{t}(\xi)$ and $\hat{Z}^{(0)}(\xi)$ are homologous as preimages of $T(\xi)\subset R(\xi,\eta)$ under $\tilde{v}^{t}$.

Now it suffices to show that $\hat{Z}^{(0)}(\xi)$ represents $h(\mathcal{E})$ in $H_r(U(\mathcal{E}),\mathbb{Z})$. To see this, we need to dive into the proof of Proposition \ref{smooth-compact-cycle} in \cite{szenes2004toric}. It turns out that the proof relies on a deformation of $\hat{Z}^{(0)}(\xi)$ to $Z(\mathbf{\epsilon})$ in $U(T_X)$ via $p(t,\cdot)^{-1}()$, where 
\begin{equation}
p_j^F(s,u) = \prod_{i=1}^n \alpha_i^F(s,u)^{\langle \alpha_i, \lambda_j\rangle},
\end{equation}
defined above Lemma 5.12 of \cite{szenes2004toric}.

Proposition 5.15 of \cite{szenes2004toric} shows that $Z_{[0,1]}$ is compact, so the image in $U(T_X)$ is bounded.
So it suffices to show that $p(t,\cdot)^{-1}() \subset \bigcap_{t\in\Delta_T}U(\mathcal{E}_t)$.
But this follows from the same argument as the previous one we use to show that $\tilde{v}^{t_1}((\tilde{v}^{t_2})^{-1}(R(\xi,\eta))) \neq 0$.
$\qed$

\begin{lem}\label{hc-integral}
(a)
$\int_{Z_\phi}\Lambda = \langle \sigma_I\rangle^{quantum}$.

(b) $\int_{Z_S}\Lambda = 0$, when $S\neq \phi$.
\end{lem}
Proof: By Corollary \ref{z-represents-h_c}, $Z_S$ represents $h(\mathcal{E})\in H_r(U(\mathcal{E}),\mathbb{Z})$. 

(a) For $Z_\phi$, we have $|q_j| < |\tilde{v}_j(u)|$. By Theorem \ref{integral-formula},
\begin{equation}\label{correlator=z_phi}
\langle\sigma_I\rangle^{quantum} =\int_{Z_\phi} \Lambda.
\end{equation}
(b) For $Z_S$, $S\neq \phi$, without loss of generality assume the index $1\in S$. Then for $u\in Z_S$ we have $|\tilde{v}_1(u)|= \exp(-\langle \xi,\beta_j\rangle )< |q_1|$. Hence $\Lambda$ can be defined on $C=\{u\in U;|\tilde{v}_1(u)| \leq e^{-\langle \xi,\beta_1\rangle}, |\tilde{v}_j(u)| = e^{-\langle \xi,\beta_j\rangle},j\geq2.\}.$ Since $Z_S = \partial C$, this leads to 
\begin{equation}
\int_{Z_S} \Lambda = \int_C d\Lambda =0.
\end{equation}
\qed

\subsection{The proof of the Main Result}\label{proof}

It is easy to show\footnote{See Proposition 6.3 of \cite{szenes2004toric}.} that 
\begin{equation*}
\displaystyle\sum_{S\subset\{1,2,...,n\}} (-1)^{|S|}T_S{\rm \ is\ homologous\ to\ } T_\delta(q)
\end{equation*}
in the open set $\{y\in  (\mathbb{C}^*)^r; y_j\neq q_j, {\rm  for\ } j=1,...,r\}$. The properness of $\tilde{v}$  (Lemma \ref{properness}) then implies that $\sum (-1)^{|S|} Z_S$ is homologous to $Z_\delta(q)$ in $U(\mathcal{E})\cap U(\beta,q)$. So
\begin{equation}
\sum_{S\subset\{1,2,...,n\}}(-1)^{|S|}\int_{Z_S}\Lambda = \int_{Z_\delta(q)}\Lambda=\int_{h_q}\Lambda.
\end{equation}

Lemma \ref{hc-integral} (b) then implies that 
\begin{equation}
\int_{Z_\phi}\Lambda =\int_{h_q}\Lambda.
\end{equation}
This together with Lemma \ref{hc-integral} (a) proves (\ref{int-hq}). Hence the Main Result is proved.

\textbf{Remark:} It is worth pointing out that
in the $(2, 2)$ case, the summation formula is further explained as a toric residue of the dual toric variety. This gives the $(2, 2)$ formula the meaning of mirror symmetry. In the $(0, 2)$ case, an explanation of the right hand side of this flavour is still lacking. In future work, we hope to describe a set of dual data and explain the right hand side as a ``$(0, 2)$ toric residue" of the dual data, making the formula into a $(0, 2)$ mirror symmetry statement.

\providecommand{\bysame}{\leavevmode\hbox to3em{\hrulefill}\thinspace}
\providecommand{\MR}{\relax\ifhmode\unskip\space\fi MR }
\providecommand{\MRhref}[2]{%
  \href{http://www.ams.org/mathscinet-getitem?mr=#1}{#2}
}
\providecommand{\href}[2]{#2}


\begin{thebibliography}{10}
\bibliographystyle{amsplain}

\bibitem{batyrev1993quantum}
Victor~V Batyrev, \emph{Quantum cohomology rings of toric manifolds}, arXiv
  preprint alg-geom/9310004 (1993).

\bibitem{batyrev2002toric}
Victor~V Batyrev and Evgeny~N Materov, \emph{Toric residues and mirror
  symmetry}, Mosc. Math. J \textbf{2} (2002), no.~3, 435--475.

\bibitem{bernstein1975number}
David~N Bernshtein, \emph{The number of roots of a system of equations},
  Functional Analysis and its Applications \textbf{9} (1975), no.~3, 183--185.

\bibitem{borisov2005higher}
Lev~A Borisov, \emph{Higher {Stanley-Reisner} rings and toric residues},
  Compositio Mathematica \textbf{141} (2005), no.~1, 161--174.

\bibitem{cox2011toric}
David~A. Cox, John~B. Little, and Henry~K. Schenck, \emph{Toric varieties},
  American Mathematical Soc., 2011.

\bibitem{donagi2014mathematical}
Ron Donagi, Josh Guffin, Sheldon Katz, Eric Sharpe, et~al., \emph{A
  mathematical theory of quantum sheaf cohomology}, Asian Journal of
  Mathematics \textbf{18} (2014), no.~3, 387--418.

\bibitem{donagi2015global-aspects}
Ron Donagi, Zhentao Lu, and Ilarion~V. Melnikov, \emph{Global aspects of (0,2)
  moduli space: Toric varieties and tangent bundles}, Communications in
  Mathematical Physics \textbf{338} (2015), no.~3, 1197--1232 (English).

\bibitem{guffin2010deformed}
Josh Guffin and Sheldon Katz, \emph{Deformed quantum cohomology and (0, 2)
  mirror symmetry}, Journal of High Energy Physics \textbf{2010} (2010), no.~8,
  1--27.

\bibitem{katz2006notes}
Sheldon Katz and Eric Sharpe, \emph{Notes on certain (0, 2) correlation
  functions}, Communications in mathematical physics \textbf{262} (2006),
  no.~3, 611--644.

\bibitem{mcorist2011revival}
Jock McOrist, \emph{The revival of (0, 2) sigma models}, International Journal
  of Modern Physics A \textbf{26} (2011), no.~01, 1--41.

\bibitem{mcorist2008half}
Jock McOrist and Ilarion~V Melnikov, \emph{Half-twisted correlators from the
  coulomb branch}, Journal of High Energy Physics \textbf{2008} (2008), no.~04,
  071.

\bibitem{melnikov2012recent}
Ilarion Melnikov, Savdeep Sethi, and Eric Sharpe, \emph{Recent developments in
  (0, 2) mirror symmetry}, SIGMA \textbf{8} (2012), 068.

\bibitem{morrison1995summing}
David~R Morrison and M~Ronen Plesser, \emph{Summing the instantons: Quantum
  cohomology and mirror symmetry in toric varieties}, Nuclear Physics B
  \textbf{440} (1995), no.~1, 279--354.

\bibitem{rojas1994convex}
J~Maurice Rojas, \emph{A convex geometric approach to counting the roots of a
  polynomial system}, Theoretical Computer Science \textbf{133} (1994), no.~1,
  105--140.

\bibitem{szenes2004toric}
Andr{\'a}s Szenes and Mich{\`e}le Vergne, \emph{Toric reduction and a
  conjecture of {Batyrev} and {Materov}}, Inventiones mathematicae \textbf{158}
  (2004), no.~3, 453--495.

\bibitem{witten1993phases}
Edward Witten, \emph{Phases of n = 2 theories in two dimensions}, Nuclear
  Physics B \textbf{403} (1993), no.~1, 159--222.

\end{thebibliography}
\end{document}